\documentclass[11pt]{article}

\usepackage{amsfonts,amsmath,amsthm,amssymb}
\usepackage{bm}
\usepackage{bbm}
\usepackage[colorlinks, citecolor=blue, urlcolor=blue, backref=page]{hyperref}
\usepackage{graphicx}
\usepackage{caption}
\usepackage{subcaption}
\usepackage{placeins}

\usepackage{graphicx}
\usepackage{geometry}
\usepackage[utf8]{inputenc}
\usepackage[T1]{fontenc}    

\usepackage{enumitem}
\usepackage{hyperref}
\usepackage{url}

\newtheorem{theorem}{Theorem}[section]

\newtheorem{proposition}[theorem]{Proposition}
\newtheorem{assumption}[theorem]{Assumption}

\theoremstyle{definition}

\newtheorem{remark}[theorem]{\textbf{Remark}}
\numberwithin{equation}{section}

\numberwithin{equation}{section}
\numberwithin{theorem}{section}

\begin{document}

\title{Large-population asymptotics for the maximum of diffusive particles with mean-field interaction in the noises}

\author{Nikolaos Kolliopoulos\thanks{Department of Mathematics, University of Michigan, \url{nkolliop@umich.edu}.}
\and
David Sanchez\thanks{Department of Mathematical Sciences, Carnegie Mellon University, \url{dsanche2@andrew.cmu.edu}.}
\and
Amy Xiao\thanks{Department of Mathematical Sciences, Carnegie Mellon University, \url{yxiao3@andrew.cmu.edu}.}
}
\maketitle

\begin{abstract}
We study the $N \to \infty$ limit of the normalized largest component in some systems of $N$ diffusive particles with mean-field interaction. By applying a universal time change, the interaction in noises is transferred to the drift terms, and the asymptotic behavior of the maximum becomes well-understood due to existing results in the literature. We expect that the normalized maximum in the original setting has the same limiting distribution as that of i.i.d copies of a solution to the corresponding McKean-Vlasov SDE and we present some results and numerical simulations that support this conjecture.
\end{abstract}


\section{Introduction}

We are interested in the $N \to \infty$ asymptotic behavior of the maximum $\max_{i \leq N}X^{i, N}$ of $N$ diffusive particles $X^{1,N}, \, X^{2, N}, \, \ldots, \, X^{N,N}$ with the following dynamics:
\begin{equation}\label{specificdrift}
\begin{aligned}
dX^{i,N}_t &= r\left(X_t^{i,N}, \, \int g^r(x)\mu^N_t(dx)\right)dt + \sigma\left(X_t^{i,N}, \, \int g^{\sigma}(x)\mu^N_t(dx)\right)dB^i_t, \qquad t \geq 0, \\
X_0^{i, N} &\sim \nu^0,
\end{aligned}
\end{equation}
for $i = 1, \, 2, \, \ldots, \, N$. In the above, $(B_t^{1})_{t \geq 0}, \, (B_t^{2})_{t \geq 0}, \, \ldots, \, (B_t^{N})_{t \geq 0}$ are independent standard Brownian motion, $r, \sigma: \mathbb{R}^2 \mapsto \mathbb{R}$ and $g^r, g^{\sigma}: \mathbb{R} \mapsto \mathbb{R}$ are given functions, $\nu^0$ is a fixed initial probability distribution, and $\mu_t^N$ is the empirical measure of the particle system which is defined as
\begin{align}
\mu_t^N = \frac{1}{N}\sum_{i=1}^N\delta_{X_t^{i, N}}
\end{align}
for all $t \geq 0$. The initial values $X_0^{1, N}, \, X_0^{2, N}, \, \ldots, X_0^{N, N}$ are assumed to be independent, and also independent from the Brownian motions $(B_t^i)_{t \geq 0}$ for $i = 1, \, 2, \, \ldots, \, N$. We consider $N$ independent copies $X^1, X^2, \ldots, X^N$ of the solution $X$ to the corresponding McKean-Vlasov equation
\begin{equation}\label{McKean0}
\begin{aligned}
dX_t &= r\left(X_t, \, \mathbb{E}[g^r(X_t)]\right)dt + \sigma\left(X_t, \, \mathbb{E}[g^{\sigma}(X_t)]\right)dB_t, \qquad t \geq 0, \\
X_0 &\sim \nu^0,
\end{aligned} 
\end{equation} 
for some standard Brownian motion $B$, and we define
\begin{equation*}
P_N^{\text{int}}(t, x) := \mathbb{P}\left(\frac{\displaystyle{\max_{i\leq N}X_t^{i,N}} - b_t^N}{a_t^N} \leq x\right) \quad \text{and} \quad P_N^{\text{ind}}(t, x) := \mathbb{P}\left(\frac{\displaystyle{\max_{i\leq N}X_t^{i}} - b_t^N}{a_t^N} \leq x\right) 
\end{equation*}
for two deterministic normalizing sequences $\{a_t^N\}_{N = 1}^{\infty}$ and $\{b_t^N\}_{N = 1}^{\infty}$. Our claim is that when $r$ and $\sigma$ are sufficiently well-behaved and for a nondegenerate probability distribution function $P(x)$ we have $P_N^{\text{ind}}(t, x) \to P(x)$ as $N \to \infty$ for all $x$, we have also $P_N^{\text{int}}(t, x) \to P(x)$ as $N \to \infty$ for all $x$, reducing our study to a problem of classical Extreme Value Theory~\cite{HF10}. When $\sigma$ is constant in its second argument, the above reduction is established in \cite{KMZ22} through a delicate analysis of the Girsanov transformation that connects the i.i.d particles $X^i$ with the particles $X^{i,N}$ that interact through their drifts. In this paper, we explore some very simple settings where mean-field interaction is also present in the noises, so the techniques developed in \cite{KMZ22} become inapplicable.

\section{Elimination of noise interaction via time change: a simple Gaussian setting}
Assume that $r(x, z) = r$ for some constant $r$, $\sigma(x, z) = \sigma(z)$ depends only on its second argument, and $\nu^0 = \text{N}(m_0, \sigma_0^2)$ for some $\sigma_0^2 > 0$. Then, \eqref{specificdrift} becomes
\begin{equation}\label{random variable}
    X^{i,N}_t = X^{i,N}_0 + rt + \int^{t}_0 \sigma\left(\int g^{\sigma}(x)\mu^N_s(dx)\right)dB^i_s, \qquad t \geq 0
\end{equation}
while the corresponding McKean-Vlasov SDE acquires the form
\begin{equation}\label{McKean1}
    X_t = X_0 + rt + \int^{t}_0 \sigma\left(\mathbb{E}\left[g^{\sigma}(X_s)\right]\right)dB_s, \qquad t \geq 0
\end{equation}
It is clear now that $X_t \sim \text{N}(m_t, \sigma_t^2)$ where $m_t = m_0 + rt$ and $\sigma_t$ satisfies 
\begin{equation}\label{var1}
\sigma_t^2 = \sigma_0^2 + \int_0^t\sigma^2\left(\mathbb{E}\left[g^{\sigma}(X_s)\right]\right)ds
\end{equation}
Under this simple setup, by \cite[Example 1.1.7]{HF10} we have that
\begin{equation}
P_N^{\text{ind}}(t, x) = \mathbb{P}\left(\frac{\displaystyle{\max_{i\leq N}X_t^{i}} - b_t^N}{a_t^N} \leq x\right) \to e^{-e^{-x}} \quad \text{as} \quad N \to \infty,
\end{equation}
i.e we have weak convergence of the normalized maximum to a standard Gumbel law, where the normalizing sequences are given by 
\begin{equation}\label{bseq}
{b^N_t} = \sigma_t \sqrt{2\log(N) - \log(\log(N)) -\log(4\pi)} + m_t
\end{equation}
and 
\begin{equation}\label{aseq}
{a^N_t} = \frac{\sigma_t}{\sqrt{2\log(N) - \log(\log(N)) -\log(4\pi)}}.
\end{equation}
To deduce a similar convergence for the normalized maximum of the interacting particles $X_t^{i,N}$ for $i = 1, \, 2, \, \ldots, \, N$, we must impose a few more conditions:
\begin{assumption}\label{ass1}
The following two conditions are in force:
\begin{enumerate}
    \item The function $x \mapsto g^{\sigma}(x)$ is Lipschitz-continuous
    \item The function $z \mapsto \sigma(z)$ is twice continuously differentiable, and there exist $m^{\sigma}, M^{\sigma} > 0$ such that $m^{\sigma} \leq \sigma(z) \leq M^{\sigma}$ and $|\sigma'(z)| + |\sigma''(z)| \leq M^{\sigma}$ for all $z \in \mathbb{R}$.
    \item There exist continuous functions $K, \tilde{K}: [0, +\infty) \mapsto \mathbb{R}_{+}$ such that for all $p \in \mathbb{N}$, $t \geq 0$, and $N \in \mathbb{N}$ it holds that
    \begin{align*}
    \mathbb{E}\left[g^{\sigma}(X_t)^{2p}\right] \leq p!K(t)^p, \qquad \mathbb{E}\left[ \left(\frac{1}{N}\sum_{i=1}^{N}g^{\sigma}(X_t^i) - \mathbb{E}\left[ g^{\sigma}(X_t)\right]\right)^{2p}\right] &\leq \frac{1}{N^p}p!K(t)^p, \\
    \text{and} \qquad \mathbb{E}\left[ \left(\frac{1}{N}\sum_{i=1}^{N}g^{\sigma}\left(X_t^{i,N}\right) - \mathbb{E}\left[ g^{\sigma}(X_t)\right]\right)^{2}\right] &\leq \frac{1}{N}\tilde{K}(t)
    \end{align*}
\end{enumerate}
\end{assumption}
\begin{remark}
The second of the three moment bounds in (iii) above is very natural and consistent with the central limit theorem, while it leads to the third bound when it is combined with the triangle inequality for the $\ell^2$ norm and with the standard estimate $\mathbb{E}[|g^{\sigma}(X_t^i) - g^{\sigma}(X_t^{i,N})|^{2}] \leq CN^{-1}$. The latter holds since $g^{\sigma}$ is Lipschitz-continuous, see e.g \cite[Theorem 2.3]{Méléard1996}.
\end{remark}
\noindent The main contribution of this paper is to prove the following result
\begin{theorem}\label{maintheorem1} 
Consider the transformation $t \mapsto \tau_N(t)$ with $$\tau_N(t) := \int^t_0  \sigma^2\left(\int g^{\sigma}(x)\mu^N_s(dx)\right)ds.$$
Under Assumption~\ref{ass1} we have that $\tau_N(t) \to \tau(t) := \int^t_0  \sigma^2\left(\mathbb{E}\left[g^{\sigma}(X_s)\right]\right)ds$ as $N \to \infty$, for any $t \geq 0$.
Moreover, for each fixed $x \in \mathbb{R}$ and $T > 0$, the function
\begin{equation}
t \mapsto \tilde{P}_N^{\text{int}}(t, x) := \mathbb{P}\left(\frac{\displaystyle{\max_{i\leq N}X_t^{i,N}} - b_{\tau^{-1}\left(\tau_N(t)\right)}^N}{a_{\tau^{-1}\left(\tau_N(t)\right)}^N} \leq x\right)
\end{equation}
converges weakly in $L^2([0, T])$ to $P(x) = e^{-e^{-x}}$ as $N \to \infty$. The latter means that for any differentiable function $\phi: [0, \, T] \mapsto \mathbb{R}$ we have the convergence
\begin{equation}\label{weak}
\int_0^T\tilde{P}_N^{\text{int}}(t, x)\phi(t)dt \to e^{-e^{-x}}\int_0^T\phi(t)dt \quad \text{as} \,\, N \to \infty.
\end{equation}
\end{theorem}

\noindent The above result has two drawbacks:
\begin{enumerate}
    \item The weak convergence alone is not enough to obtain $\tilde{P}_N^{\text{int}}(t, x) \to e^{-e^{-x}}$ for fixed $t$ and $x$, which is required to deduce that the normalized maximum at any $t$ converges weakly to a standard Gumbel distribution.

    \item The normalizing sequences that are used are $\left\{a_{\tau^{-1}\left(\tau_N(t)\right)}^N\right\}_{N=1}^{\infty}$ and $\left\{b_{\tau^{-1}\left(\tau_N(t)\right)}^N\right\}_{N=1}^{\infty}$, which are non-deterministic.
\end{enumerate}
However, the next proposition shows that if the weak convergence in $L^2([0, T])$ can be upgraded to strong convergence in $L^2([0, T])$, we will also be able to replace the stochastic normalizing sequences by the natural deterministic candidates $\left\{a_{t}^N\right\}_{N=1}^{\infty}$ and $\left\{b_{t}^N\right\}_{N=1}^{\infty}$. As a consequence, the desired weak convergence of the normalized maxima of $X_t^{i,N}$ reduces to showing that the sequence of functions $\{\tilde{P}_N^{\text{int}}(\cdot, x)\}_{N=1}^{\infty}$ is e.g relatively compact (which requires uniform equicontinuity in $[0, T]$). We will focus on this problem in a follow-up work, where we will also focus on covering a much wider class of mean-field systems. 

\begin{proposition}\label{maintheorem2}
With the notation of Theorem~\ref{maintheorem1}, assume that for some $t \geq 0$ we have the following convergence:
\begin{equation}\label{strong}
\tilde{P}_N^{\text{int}}(t, x) = \mathbb{P}\left(\frac{\displaystyle{\max_{i \leq N}X_t^{i,N}} - b_{\tau^{-1}\left(\tau_N(t)\right)}^N}{a_{\tau^{-1}\left(\tau_N(t)\right)}^N} \leq x\right) \to e^{-e^{-x}} \quad \text{as} \,\, N \to \infty,
\end{equation}
for all $x \in \mathbb{R}$. Then, we have also the convergence:
\begin{equation}
P_N^{\text{int}}(t, x) = \mathbb{P}\left(\frac{\displaystyle{\max_{i \leq N}X_t^{i,N}} - b_{t}^N}{a_{t}^N} \leq x\right) \to e^{-e^{-x}} \quad \text{as} \,\, N \to \infty,
\end{equation}
for all $x \in \mathbb{R}$.
\end{proposition}
\noindent We proceed now to the proofs of the above results.
\begin{proof}[Proof of Theorem~\ref{maintheorem1}]
For each i, we define the process
\begin{equation*}
W_t^i = \int^{\tau^{-1}_N(t)}_0 \sigma\left(\int g^{\sigma}(x)\mu^N_s(dx)\right) dB_s^i
\end{equation*} 
so that $(W^1)_{t \geq 0}, (W^2)_{t \geq 0}, \ldots, (W^N)_{t \geq 0}$ are $N$ independent standard Brownian motions. Hence we can rewrite the initial formulation as 
\begin{equation}\label{timechanged1}
X^{i,N}_t = X^{i,N}_0 + rt + W^i_{\tau_{N}(t)}
\end{equation}
Plugging $t \rightarrow \tau^{-1}_N(t)$ we get:
\begin{align}\label{timechangedX}
    X^{i,N}_{\tau^{-1}_N(t)} &= X^{i,N}_0 + r\tau^{-1}_N(t) + W^i_t, \qquad t \geq 0 
\end{align}
Note that $\tau_N$ is a monotonically increasing and differentiable function with $$(\tau^{-1}_N(t))' = \frac{1}{\tau'_N(\tau^{-1}_N(t))} = \frac{1}{\sigma^2\left(\int g^{\sigma}(x)\mu^N_{\tau^{-1}_N(t)}(dx)\right)} = \frac{1}{\sigma^2\left(\frac{1}{N} \displaystyle{\sum_{i=1}^N}g\left(X^{i,N}_{\tau^{-1}_N(t)}\right)\right)}.$$
Thus, setting $Y^{i,N}_t := X^{i,N}_{\tau^{-1}_N(t)}$, equation \eqref{timechangedX} becomes
\begin{equation}\label{Yparticle}
Y^{i,N}_{t} = X^{i,N}_0 + r \int^{t}_0 \frac{1}{\sigma^2\left(\int g^{\sigma}(Y)\mu^{Y,N}_s(dy)\right)}ds + W^i_t, \quad \quad t \geq 0,  
\end{equation}
where $\mu^{Y,N} = \displaystyle{\sum_{i=1}^N}g^{\sigma}\left(Y^{i,N}_{t}\right)$ is the empirical measure of the timechanged particles. The time-changed mean-field system ~\eqref{Yparticle}
has corresponding McKean-Vlasov equation 
\begin{equation}\label{YMckean}
    Y_t = X_0 + r\int_0^t\frac{1}{\sigma^2\left(\mathbb{E}[g^{\sigma}\left(Y_s\right)]\right)}ds + W_t, \quad \quad t \geq 0,
\end{equation}
with $X_0 \sim \text{N}(m_0, \sigma^2_0)$ and $W$ being a standard Brownian motion. This yields
\begin{equation}\label{VarSDE}
\mathbb{E}[Y_t] = m_0 + r\int_0^t\frac{1}{\sigma^2\left(\mathbb{E}[g^{\sigma}\left(Y_s\right)]\right)}ds \quad \text{and} \quad \text{Var}(Y_t) = \sigma_0^2 + t
\end{equation}

\noindent Observe now that the particles in the mean-field system \ref{Yparticle} interact only through the drifts, and they have a normal distribution which belongs to the domain of attraction of the Gumbel extreme value distribution. Moreover, the conditions of \cite[Theorem 2.4]{KMZ22} are applicable due to Assumption~\ref{ass1} and we can argue as in \cite[Example 3.1]{KMZ22} to obtain
\begin{equation}\label{GumbelY}
\mathbb{P}\left(\frac{\displaystyle{\max_{i \leq N}Y^{i,N}_t} - {\tilde{b}^N_t}}{\tilde{a}^N_t} \leq x\right) \to e^{-e^{-x}} \quad \text{as} \quad N \to \infty,
\end{equation}
for all $t \geq 0$, where the normalizing sequences are given by $$\tilde{b}^N_t =  (\sigma_0^2 + t)\sqrt{2\log(N) - \log(\log(N)) -\log(4\pi)} + m_0 + r\int_0^t\frac{1}{\sigma^2\left(\mathbb{E}[g^{\sigma}\left(Y_s\right)]\right)}ds$$ and $$\tilde{a}^N_t = \frac{\sigma_0^2 + t}{\sqrt{2\log(N) - \log(\log(N)) -\log(4\pi)}}.$$ 
By (ii) and the third bound in (iii) of Assumption~\ref{ass1} and the dominated convergence theorem, we can easily obtain that $\tau_N(t)$ converges for every $t \geq 0$ to the deterministic time-changing function $\tau(t) = \int_0^t\sigma^2(\mathbb{E}[g^{\sigma}(X_s)])ds$, where $(X_t)_{t \geq 0}$ solves \eqref{McKean1}. Thus, the natural coupling for the particles $X_t^{i, N} = Y_{\tau_N(t)}^{i, N}$ are the particles $Y_{\tau(t)}^{i, N}$. We will show now that
\begin{equation}\label{weakstochseq}
\int_0^T\mathbb{P}\left(\frac{\displaystyle{\max_{i \leq N}X_t^{i,N}} - \tilde{b}_{\tau_N(t)}^N}{\tilde{a}_{\tau_N(t)}^N} \leq x\right)\phi(t)dt - \int_0^T\mathbb{P}\left(\frac{\displaystyle{\max_{i \leq N}Y_{\tau(t)}^{i,N}} - \tilde{b}_{\tau(t)}^N}{\tilde{a}_{\tau(t)}^N} \leq x\right)\phi(t)dt \to 0
\end{equation}
as $N \to \infty$. Approximating $\phi$ with a linear combination of simple indicator functions, the above reduces to showing
\begin{equation}\label{difference}
\mathbb{E}\left[\int^T_0 \prod^{N}_{i=1} \mathbbm{1}_{\{ X^{i,N}_t \leq \tilde{a}^N_{\tau_N(t)} x + \tilde{b}^N_{\tau_N(t)} \}}dt\right] - \mathbb{E}\left[\int^T_0 \prod^{N}_{i=1} \mathbbm{1}_{\{ Y^{i,N}_{\tau(t)} \leq \tilde{a}^N_{\tau(t)} x + \tilde{b}^N_{\tau(t)}\}}dt\right] \to 0 
\end{equation}
for arbitrary $T > 0$. We perform the change of variable $t \rightarrow \tau_N^{-1}(\tau(t))$ on the first integral to write it as 
\begin{align*}
&\mathbb{E}\left[\int^{\tau^{-1}(\tau_N(T))}_0 \prod^{N}_{i=1} \mathbbm{1}_{\{ Y^{i,N}_{\tau(t)} \leq \tilde{a}^N_{\tau(t)} x + \tilde{b}^N_{\tau(t)}\}}\frac{\tau'(t)}{\sigma^2\left(\int g^{\sigma}(x)\mu^N_{\tau^{-1}_N(\tau(t))}(dx)\right)}dt\right] \\
& \qquad \qquad \qquad = \mathbb{E}\left[\int^{\tau^{-1}(\tau_N(T))}_0 \prod^{N}_{i=1} \mathbbm{1}_{\{ Y^{i,N}_{\tau(t)} \leq \tilde{a}^N_{\tau(t)} x + \tilde{b}^N_{\tau(t)}\}}\frac{\sigma^2(\mathbb{E}[g^{\sigma}(X_t)])}{\sigma^2\left(\int g^{\sigma}(Y)\mu^{Y,N}_{\tau(t)}(dy)\right)}dt\right].
\end{align*}
Then, the absolute value of \eqref{difference} can be written as
\begin{align*}
&\Bigg|\mathbb{E}\left[\int^{\tau^{-1}(\tau_N(T))}_T \prod^{N}_{i=1} \mathbbm{1}_{\{ Y^{i,N}_{\tau(t)} \leq \tilde{a}^N_{\tau(t)} x + \tilde{b}^N_{\tau(t)}\}}\frac{\sigma^2(\mathbb{E}[g^{\sigma}(X_t)])}{\sigma^2\left(\int g^{\sigma}(Y)\mu^{Y,N}_{\tau(t)}(dy)\right)}dt\right] \\
&\qquad \qquad \qquad + \mathbb{E}\left[\int^{T}_0 \prod^{N}_{i=1} \mathbbm{1}_{\{ Y^{i,N}_{\tau(t)} \leq \tilde{a}^N_{\tau(t)} x + \tilde{b}^N_{\tau(t)}\}}\left(\frac{\sigma^2(\mathbb{E}[g^{\sigma}(X_t)])}{\sigma^2\left(\int g^{\sigma}(Y)\mu^{Y,N}_{\tau(t)}(dy)\right)} - 1\right)dt\right]\Bigg| \\
& \quad \qquad \leq C\left(\mathbb{E}[|\tau^{-1}(\tau_N(T)) - T|] + \mathbb{E}\left[\int^{T}_0\left(\frac{\sigma^2(\mathbb{E}[g^{\sigma}(X_t)])}{\sigma^2\left(\int g^{\sigma}(Y)\mu^{Y,N}_{\tau(t)}(dy)\right)} - 1\right)dt\right]\right) 
\end{align*}
with $C$ depending on the positive upper and lower bounds of the function $\sigma$. By the continuity of $\tau$ and thus of $\tau^{-1}$ and the convergence $\tau_N(t) \to \tau(t)$, we get the almost sure convergence $\tau^{-1}(\tau_N(T)) \to T$ as $N \to \infty$. Thus, by the boundedness of the function $\sigma$ and the dominated convergence theorem, the bound in the previous display goes to
\begin{equation*}
C \mathbb{E}\left[\int^{T}_0\left(\frac{\sigma^2(\mathbb{E}[g^{\sigma}(X_t)])}{\sigma^2(\mathbb{E}[g^{\sigma}(Y_{\tau(t)})])} - 1\right)dt\right] = 0
\end{equation*}
since we can plug $t \rightarrow \tau^{-1}(t)$ on \eqref{McKean1} to find that $X_{\tau^{-1}(t)}$ satisfies \eqref{YMckean}, for which weak uniqueness holds by \cite[Proposition 1]{Chaintron_2022} and (i) of Assumption~\ref{ass1}, so $X_{t}$ must have the same law as $Y_{\tau(t)}$. This finishes the proof of \eqref{weakstochseq}.

Observe now that \eqref{var1} gives $\sigma_t^2 = \sigma_0^2 + \tau(t) \Rightarrow \sigma_0^2 + t = \sigma_{\tau^{-1}(t)}^2$, which implies that $\tilde{a}_{t}^N = a_{\tau^{-1}(t)}^N$ and $\tilde{b}_{t}^N = b_{\tau^{-1}(t)}^N$. Hence we have that
\begin{equation*}
\mathbb{P}\left(\frac{\displaystyle{\max_{i \leq N}X_t^{i,N}} - \tilde{b}_{\tau_N(t)}^N}{\tilde{a}_{\tau_N(t)}^N} \leq x\right) = \tilde{P}_N^{\text{int}}(t, x).
\end{equation*}
Plugging the above in \eqref{weakstochseq} and recalling also \eqref{GumbelY}, we obtain \eqref{weak} and our proof is complete.
\end{proof}

\begin{proof}[Proof of Proposition~\ref{maintheorem2}]
Suppose that $\tilde{P}_N^{\text{int}}(t, x) \to e^{-e^{-x}}$ as $N \to \infty$ for all $x \in \mathbb{R}$, for a fixed $t \geq 0$. This is equivalent to
$$\frac{\displaystyle{\max_{i \leq N}X_t^{i,N}} - b_{\tau^{-1}\left(\tau_N(t)\right)}^N}{a_{\tau^{-1}\left(\tau_N(t)\right)}^N}$$
converging weakly to the standard Gumbel law as $N \to \infty$. Now we have
\begin{equation*}
\frac{a_{\tau^{-1}(\tau_N(t))}^N}{a_{t}^N} = \frac{\sigma_{\tau^{-1}(\tau_N(t))}}{\sigma_{t}} = \sqrt{\frac{\sigma_0^2 + \tau_N(t)}{\sigma_0^2 + \tau(t)}} \to 1 \quad \text{as} \quad N \to \infty
\end{equation*}
Thus, we can deduce that
\begin{equation*}
\frac{\displaystyle{\max_{i \leq N}X_t^{i,N}} - b_{\tau^{-1}(\tau_N(t))}^N}{a_{t}^N} = \frac{a_{\tau^{-1}(\tau_N(t))}^N}{a_{t}^N}\frac{\displaystyle{\max_{i \leq N}X_t^{i,N}} - b_{\tau^{-1}(\tau_N(t))}^N}{a_{\tau^{-1}(\tau_N(t))}^N}
\end{equation*}
converges also weakly to the standard Gumbel law. Then we can write
\begin{equation*}
\frac{\displaystyle{\max_{i \leq N}X_t^{i,N}} - b_{t}^N}{a_{t}^N} = \frac{\displaystyle{\max_{i \leq N}X_t^{i,N}} - b_{\tau^{-1}(\tau_N(t))}^N}{a_{t}^N} + \frac{b_{\tau^{-1}(\tau_N(t))}^N - b_{t}^N}{a_{t}^N}
\end{equation*}
which will also converge weakly to the standard Gumbel law, since we have
\begin{align*}
\left|\frac{b_{\tau^{-1}(\tau_N(t))}^N - b_{t}^N}{a_{t}^N}\right| &= \left|\frac{\sigma_{\tau^{-1}(\tau_N(t))}^2 - \sigma_{t}^2}{\sigma_{t}(\sigma_{\tau^{-1}(\tau_N(t))} + \sigma_{t})}\left(2\log(N) - \log(\log(N)) -\log(4\pi)\right)\right| \\
&\leq C\log(N)\left|\tau_N(t) - \tau(t)\right| \\
&= C\int_0^t\log(N)\left|\sigma^2\left(\int g^{\sigma}(x)\mu^N_s(dx)\right) - \sigma^2\left(\mathbb{E}\left[g^{\sigma}(X_s)\right]\right)\right|ds
\end{align*}
which goes to $0$ by (ii) and the third moment bound in (iii) of Assumption~\ref{ass1}. The last weak convergence is precisely $P_N^{\text{int}}(t, x) \to e^{-e^{-x}}$ as $N \to \infty$ for all $x \in \mathbb{R}$
\end{proof}

\section{A model for a network of banks: numerical simulations}
Systems of the form \eqref{specificdrift} have been used for modeling the monetary reserves of banks that interact through borrowing, in which case the drift term has the form $$r\left(X_t^{i,N}, \, \int g^r(x)\mu^N_t(dx)\right) = \kappa\left(X^{i,N}_t - \frac{1}{N}\sum^{N}_{j=1}X^{j,N}_t\right)$$ for $\kappa < 0$, which captures the tendency of a bank to borrow money when its reserves are below the average wealth of all the banks, and to lend money otherwise. We refer to \cite{BC2015} for a system with additional control terms in the drifts which are chosen by the banks to minimize an expected cost (a Mean-Field game), and to \cite{CFS15} for a setting where each particle has an idiosyncratic CIR-like volatility. Studying the large-$N$ asymptotics of $\max_{i \leq N}X^{i, N}$ as $N \to \infty$ in all those settings could lead to a better understanding of the behaviour of the richest banks in large interbanking networks, and also the poorest ones since the dynamics are invariant upon a multiplication with $-1$ and we clearly have $\min_{i \leq N}X^{i, N} = -\max_{i \leq N}(-X^{i, N})$. Our aim in the future is to adapt the whole theory of extreme values to the framework of mean-field systems, developing also statistical methods (see \cite[Chapters 3-4]{HF10}) which may be helpful for estimating credit risk under the above models.

We would like to check whether the ideal result $P_N^{\text{int}}(t, x) \to e^{-e^{-x}}$ as $N \to \infty$ for all $x \in \mathbb{R}$ holds under the following model for a network of banks

\begin{equation}\label{new rv}
X_t^{i, N} = X_0^{i, N} + \int_0^t\left(X_s^{i, N} - \frac{1}{N}\sum_{j = 1}^NX_s^{j, N}\right)ds + \int_0^t\sqrt{\frac{1}{N}\sum_{j = 1}^N\left(X_s^{j, N}\right)^2}dB_s^i.
\end{equation}
with the initial values being i.i.d random variables with a $\text{N}(0, 1)$ distribution. In this case we can easily find that \eqref{McKean0} acquires the form $dX_t = (X_t - \mathbb{E}[X_t])dt + \sqrt{\text{Var}(X_t)}dB_t$, where we can take expectations to get $\mathbb{E}[X_t] = 0$, and then apply Ito's formula on $X_t^2$ and take expectations to obtain $$\text{Var}(X_t) = 1 + 3\int_0^t\text{Var}(X_s)ds$$ and thus $\text{Var}(X_t) = e^{3t}$. Hence, we can deduce that $X_1 \sim \text{N}(0, e^{3})$. Obtaining even partial results like Theorem~\ref{maintheorem1} and Proposition~\ref{maintheorem2} is difficult for this setting due to the insufficient boundedness of the SDE coefficients and their derivatives. For this reason, we will restrict to a numerical analysis. We used the Euler-Maruyama method with step $\Delta t=0.0001$ on \eqref{new rv} to simulate approximations $\tilde{X}_1^{i,N}$ for the variables $X_1^{i,N}$ for $i \in \{1, 2, \ldots, N\}$, and we evaluated the approximate normalized maximum
\begin{equation}\label{maxapprox}
M^N_t = \frac{\displaystyle{\max_{i\leq N}\tilde{X}_1^{i, N}} - b_1^{N}}{a_1^{N}} = \frac{\displaystyle{\max_{i\leq N}X_1^{i, N}} - b_1^{N}}{a_1^{N}} + \text{Error}
\end{equation}
with $a_1^{N}$ and $b_1^{N}$ given by \eqref{aseq} and \eqref{bseq} for $t = 1$, where $m_1 = 0$ and $\sigma_1 = e^{\frac{3}{2}}$. We repeated this process $1000$ times to simulate $1000$ independent normalized maxima $M_1, M_2, \ldots, M_{1000}$ and we drew a histogram of the values $U_i = F(M_i)$ for $F(x) = e^{-e^{-x}}$, ~using~ intervals of length $0.1$. By \cite[Theorem 10.2.2]{KP92}, the error in \eqref{maxapprox} for $N \leq 200$ is at most of order $\sqrt{\Delta t} \times \sqrt{2\log(200)} \approx 0.03$. If $p_i$ is the probability that $U_t^N = F(M_t^N)$ falls in the $[\frac{i-1}{10}, \frac{i}{10}]$, the error due to the law of large numbers is of order $\sqrt{0.001} \times \sqrt{p_i(1 - p_i)} \leq 0.015$. Thus, the histogram of the $U_i$ is a good approximation of the distribution of $U_t^N$, which should be the uniform distribution if the desired result is true and $N$ is sufficiently large. We first conducted the above analysis twice for $N = 100, 141, 173$, and $200$, obtaining histograms that resemble a uniform distribution in $[0, 0.9]$, with a bigger error in $[0.9, 1]$:   

\FloatBarrier
\begin{figure}[h]
    \centering
    \includegraphics[width=15cm]{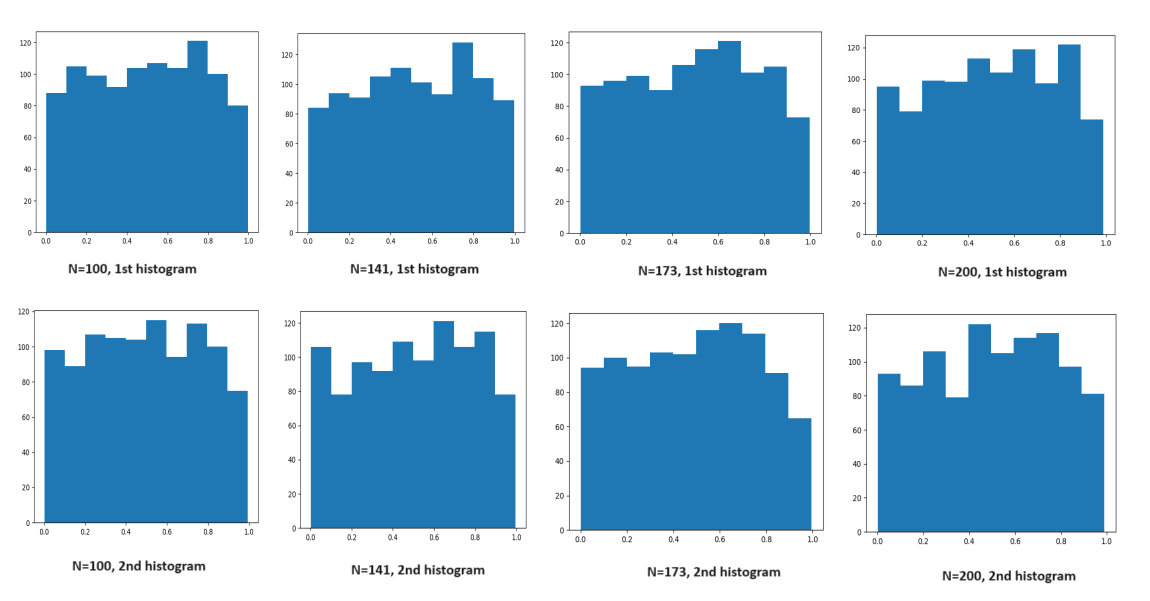}
\end{figure}
\FloatBarrier

\noindent Considering the error due to the law of large numbers in the above graphs, it seems that the distribution of $U_t^N$ tends to come close to the uniform distribution on $[0, 1]$ for $N > 100$. This indicates that the desired result is probably true, with a quite fast rate of convergence. 

Next, we repeat the above process with the variables $\tilde{X}_1^{i,N}$ replaced with the i.i.d McKean-Vlasov particles $X_1^{i} \sim N(0, e^3)$. Since normal random variables can be simulated directly, the complexity is much smaller and allows us to take much larger values of $N$ and a larger number of simulated maxima which leads to more accurate pictures. The following histograms are derived for $N=100, 1000, 10000$ and $100000$, with $10000$ simulated maxima:

\FloatBarrier
\begin{figure}[h]
    \centering
    \includegraphics[width=15cm]{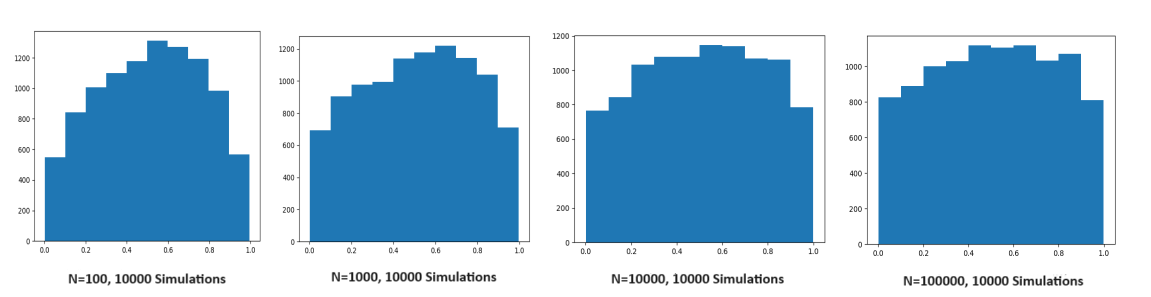}
\end{figure}
\FloatBarrier

\noindent The above shows convergence to a uniform distribution, which is known to hold by standard Extreme Value Theory, but good approximations require $N > 10000$ if not $N > 100000$. This is consistent with the result of \cite{Nair}, where a slow logarithmic rate of convergence is derived for normalized maxima of independent Gaussian variables. However, the interesting finding is that the distribution of the normalized maximum of the interacting particles seems to converge much faster than that of the i.i.d McKean-Vlasov particles, with accurate approximations obtainable for $N \leq 200$, at least for $F(M_t^N) < 0.9$. Finally, it seems that it is the introduction of Mean-Field interaction in the noise terms that leads to faster convergence. Indeed, repeating the first numerical experiment but with only the volatility term $\sqrt{\frac{1}{N}\sum_{j = 1}^N\left(X_s^{j, N}\right)^2}$ replaced with its limit $e^{\frac{3s}{2}}$, that is
\begin{equation}\label{newnew rv}
X_t^{i, N} = X_0^{i, N} + \int_0^t\left(X_s^{i, N} - \frac{1}{N}\sum_{j = 1}^NX_s^{j, N}\right)ds + \int_0^te^{\frac{3s}{2}}dB_s^i,
\end{equation}
the following histograms are generated for $N=150$ with $1000$ simulated maxima, where the distribution of $U_t^N = F(M_t^N)$ seems to deviate significantly from the uniform in $[0, 1]$:

\FloatBarrier
\begin{figure}[h]
    \centering
    \includegraphics[width=12.5cm]{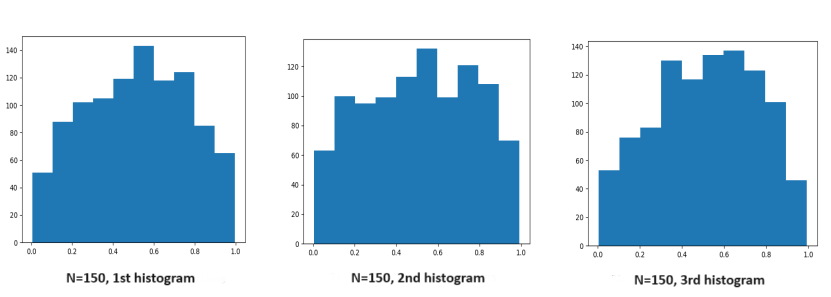}
\end{figure}
\FloatBarrier

\bibliography{references}

\begin{thebibliography}{1}

\bibitem{BC2015}
L.~Bo and A.~Capponi.
\newblock Systemic risk in interbanking networks.
\newblock {\em SIAM J. Financial Math.}, 6(1):386--424, 2015.

\bibitem{CFS15}
R.~Carmona, J.~Fouque, and L.~Sun.
\newblock Mean field games and systemic risk.
\newblock {\em Communications in Mathematical Sciences}, 13(4):911--933, 2015.

\bibitem{Chaintron_2022}
L.-P. Chaintron and A.~Diez.
\newblock Propagation of chaos: A review of models, methods and applications.
  i. models and methods.
\newblock {\em Kinetic and Related Models}, 15(6):895, 2022.

\bibitem{HF10}
L.~de~Haan and A.~Ferreira.
\newblock {\em Extreme Value Theory: An Introduction}.
\newblock Springer Series in Operations Research and Financial Engineering.
  Springer New York, 2007.

\bibitem{KP92}
P.~E. Kloeden and E.~Platen.
\newblock {\em Numerical solution of stochastic differential equations},
  volume~23 of {\em Applications of Mathematics (New York)}.
\newblock Springer-Verlag, Berlin, 1992.

\bibitem{KMZ22}
N.~Kolliopoulos, M.~Larsson, and Z.~Zhang.
\newblock Propagation of chaos for maxima of particle systems with mean-field
  drift interaction.
\newblock {\em Probability Theory and Related Fields}, 2023.

\bibitem{Méléard1996}
S.~M{\'e}l{\'e}ard.
\newblock {\em Asymptotic behaviour of some interacting particle systems;
  McKean-Vlasov and Boltzmann models}, pages 42--95.
\newblock Springer Berlin Heidelberg, Berlin, Heidelberg, 1996.

\bibitem{Nair}
K.~A. Nair.
\newblock {Asymptotic Distribution and Moments of Normal Extremes}.
\newblock {\em The Annals of Probability}, 9(1):150 -- 153, 1981.

\end{thebibliography}
\bibliographystyle{abbrv}

\end{document}